\hyphenation{para-met-rize}

\def\union{\cup}
\def\W{{\cal W}}
\def\P{{\bf P}}

\def\Q{{\bf Q}}

\def\ref#1#2#3#4#5#6{ #1: #2. #3 {\bf #4}, #5 (#6)}

\def\card{\mathop{\rm card}}

\def\Pic{\mathop{\rm Pic}}
\def\mgbar{{\overline{\cal M}_g}}
\def\mgb#1{\overline{\cal M}_{#1}}
\def\mgn{{\cal M}_{g,n}}
\def\mbar#1#2{\overline{\cal M}_{#1,#2}}
\def\mgnbar{\mbar{g}{n}}
\def\define{\smallskip{\noindent\bf Definition. }}
\def\which#1#2{\smallskip\global\advance \thmnum by 1 \edef#2{#1 \the\secnum .\the\thmnum}\noindent{\bf #2}\ \thinspace}
\def\thm#1{\which{Theorem}{#1}}
\def\lemma#1{\which{Lemma}{#1}}
\def\cor#1{\which{Corollary}{#1}}

\def\prop#1{\which{Proposition}{#1}}
\def\proof{\smallskip{\noindent\sl Proof.\ }}
\def\floor#1{\lfloor #1 \rfloor}
\def\ref#1#2#3#4#5#6{ #1: #2. #3 {\bf #4}, #5 (#6)}
\def\mapright#1{\smash{\mathop{\longrightarrow}\limits^{#1}}}
\def\mapdown#1{\Big\downarrow\rlap{$\vcenter{\hbox{$\scriptstyle#1$}}$}}

 4
\font\largebold=cmbx10 scaled \magstep 3
 2
\font\big=cmbx12
\font\sit=cmti8

\def\Q{{\bf Q}}

\def\isisom{\cong}
\def\chapter#1{\vfill\eject{\largebold\centerline{#1}\par}\vskip .2in}
\def\section#1{\smallskip\penalty -1000\global\advance \secnum by 1 
{\big\centerline{{\ifnum\the\secnum>0 \the\secnum. \fi}#1}\par}\vskip .1in \global\thmnum=0}

\overfullrule=0pt
\def\mapsto{\to}
\def\mapsby#1{\smash{\mathop{\mapsto}\limits^{#1}}}

\def\widen{\def\normalbaselines{\baselineskip20pt \lineskip3pt \lineskiplimit 3pt}}
\newcount\secnum
\newcount\thmnum
\thmnum = 0
\def\qed{ \ {\vbox{\hrule height 0pt depth .3pt\hbox{\vrule width .3pt height 6pt
\kern 3pt \vrule width .3pt}\hrule height 0pt depth .3pt}}\smallskip}

\def\title#1{{\centered{\largebold #1}\medskip}}
\def\centered{\leftskip=0pt plus 4 em \rightskip=\leftskip
\parfillskip=0 pt \spaceskip=.3333 em \xspaceskip=.5em
\pretolerance=9999 \tolerance=9999 \parindent=0pt \hyphenpenalty=5000
\exhyphenpenalty=5000}
\def\today{\ifcase\month\or January\or February\or March\or April\or May\or
June\or July\or August\or September\or October\or November\or December\fi
\ \number\day, \number\year}

\title{Relations among divisors on the moduli space of curves with marked
points}
\centerline{Adam Logan}
\smallskip
\centerline{\today}
\medskip
\section{Introduction}
Let $\mgbar$ be the coarse moduli space of stable curves of genus $g$.
Eisenbud, Harris and Mumford proved a relation between certain divisors
on $\mgbar$ (the Brill-Noether divisors, to be described below).  Calculating
their classes in $\Pic \mgbar \otimes \Q$, they succeeded in proving that
$\mgbar$ is of general type for $g > 23, g+1$ composite.  In subsequent
work, the restriction that $g+1$ be composite was removed.

In this paper, I will generalize their relation to $\mgnbar$, the moduli space
of stable curves of genus $g$ with $n$ marked points.  This does not yield new
results on the Kodaira dimension of $\mgnbar$, as the ``divisors of
Brill-Noether type'', which I introduce below, are less effective for this
purpose than certain other divisors which I studied in my dissertation [L].
(I expect to publish the other results of [L] shortly.)

The remainder of this introduction will be devoted to stating the results.
For basic facts on $\mgbar$ and $\mgnbar$, the reader is referred
to [HM] and [K] respectively.

\define Fix a nonnegative integer $g$, and let $r > 1, d > 1$ be integers such
that $g - (r+1)(g-d+r) = -1$.  (Here, the left side is the expected dimension
of the space of $g^r_d$'s on a curve of genus $g$.)  A {\sl divisor of
Brill-Noether type} on $\mgbar$ is a codimension-1 component of the locus of
curves which have an admissible $g^r_d$.

(All $g^r_d$'s on nonsingular irreducible curves are admissible;
in general, an admissible
$g^r_d$ is a $g^r_d$ on each component with certain ramification conditions
at the singular points of the curve.  Again, see [HM] for details.)

The Brill-Noether Ray Theorem of Eisenbud, Harris, and Mumford then asserts
that, for fixed $g$, all of these divisors on $\mgbar$ are linearly dependent,
and calculates their class.  We generalize their definition and theorem as
follows:

\define Fix $g$ and $n$, and let integers $r, d$ and sequences
$\{Z_1\}, \dots, \{Z_n\}$ of length $r+1$ be such
that $$g-(r+1)(g-d+r)-\sum Z_{i,j} = -1.$$
(The left-hand side is the expected dimension of the set of $g^r_d$'s
on a curve of genus $g$ with ramification sequences $\{Z_1\}, \dots,
\{Z_n\}$ at
points $p_1, \dots, p_n$.)  A {\sl divisor of Brill-Noether type} on
$\mgnbar$ is a divisorial component of the locus of curves and sets of points
such that there is a $g^r_d$ on the curve whose ramification sequences
at the $p_i$ are at least $\{Z_i\}$.

\thm{\bntthm}
For any $n > 2$, every divisor of Brill-Noether type on $\mgnbar$ is a linear
combination of pullbacks of divisors of Brill-Noether type from
$\mbar{g}{2}$.  Also, the space of divisors of Brill-Noether type on
$\mbar{g}{1}$ has dimension $2$, for any $g > 2$.\qed

\cor{\bntdim}
The dimension of the subspace of $\Pic \mgnbar \otimes \Q$ spanned by the
divisors of Brill-Noether type is $1 + n + {n \choose 2}$.\qed

{\sl Acknowledgments.}  This paper is a lightly revised version of parts of
my 1999 doctoral dissertation [L], written under the direction of Joe Harris,
to whom I am very grateful.  I would also like to thank Mira Bernstein for
teaching me how to work with moduli spaces.
Support for this work was provided by a Sloan Dissertation Fellowship.

\section{The First Case}
Throughout the paper, we work over an algebraically closed field of 
characteristic 0 (though this is surely unnecessary), and we deal only with
genus $\ge 3$.  We start by introducing some notation.  Arbarello and Cornalba
have given a basis for $\Pic \mgnbar \otimes \Q$ [AC, Thm.\ 2]:

\thm{\picgen} $\Pic_{\hbox{\sit fun}} \mgnbar$, the Picard group of the moduli
stack $\mgnbar$, is free on the following
generators: $\lambda, \delta_0, \psi_i\, (1 \le i \le n),$ and
$\delta_{i;S}\, (0 \le i \le \floor{g/2}), S \subseteq \{1,2,\dots,n\}),
\card S > 1\hbox{ if }i = 0$.\qed

Here $\lambda$ is the pullback of $\lambda$ on $\mgbar$, $\delta_0$ the
divisor corresponding to the locus of curves with a nondisconnecting node,
$\delta_{i,S}$ the divisor corresponding to the locus of curves with a
node whose removal leaves one component of genus $i$ with precisely the
marked points indexed by $S$, and $\psi_i$ is the divisor class on the
moduli stack which takes the value $-\pi_* (\sigma_i^2)$ on the family
${\cal X} \, \mapsby {\pi} \, B$ with the $\sigma_i$ as sections.

Usually we will prefer to work with $\omega_i$, the relative dualizing sheaf
of the $i\/$th projection map from $\mgnbar$ to $\mbar{g}{n-1}$.  If we
replace $\psi_i$ by $\omega_i$ in the above, the theorem remains true,
because $\omega_i = \psi_i - \sum_{i \in S} \delta_{0;S}$, as will be seen
below.

We will be pulling back divisors, so it also seems appropriate to state the
results concerning this.  Again, the answer is given in [AC, p.\ 161]:

\thm{\pullback} The pullback map is given as follows:
$$\eqalign{\pi_n^* \lambda &= \lambda,\cr
\pi_n^*\delta_0 &= \delta_0,\cr
\pi_n^* \omega_i &= \omega_i,\cr
\pi_n^* \psi_i &= \psi_i - \delta_{0;\{i,n\}}\cr
\pi_n^* \delta_{i;S} &= \delta_{i;S} + \delta_{i;S \union \{n\}},\cr}$$
except that $$\pi_1^* \delta_{g/2;\emptyset} = \delta_{g/2;\emptyset}$$ for
$n = 1$.\qed

For notational simplicity, we have only stated this for $\pi_n$, but the
action of the symmetric group ${\cal S}_n$ on $\mgnbar$ makes it easy to
describe the effects of the other $\pi$.  Also observe that, by an easy
induction starting from $\psi = \omega$ on $\mbar{g}{1}$, we get
$$\psi_i = \omega_i + \sum_{i \in S \subseteq \{1,2,\dots,n\}, S \ne \{i\}}
\delta_{0;S}$$ on $\mgnbar$, as asserted above.

We now start by dealing with the case $n = 1$.

\define Let $BN$ be the divisor class $$(g+3) \lambda -{(g+1)\over 6} \delta_0
- \sum_{i=1}^{\floor{g/2}} i(g-i) \delta_i$$ on $\mgbar$.

If $g+1$ is not prime, then $BN$ is a positive multiple of the class of an
effective divisor, namely that of the union of the codimension-$1$ components
of the locus of curves admitting a $g^r_d$, where $(r+1)(g-d+r) = g+1$ [EH, 
Thm.\ 1] (the condition on $r, d$, and $g$ ensures that there are such
components).

This will be one of the generators of the Brill-Noether space on $\mbar{g}{1}$
when $g+1$ is not prime; the other will be
the Weierstrass divisor, that is, the locus of Weierstrass points.  Its class
was computed by Cukierman [Cuk, Thm.\ 2.0.12 and following remark].  His
result is as follows:

\thm{\cukwei}
Let $g \ge 2$.  The class on $\mbar{g}{1}$ of the {\sl Weierstrass 
divisor} $\W$, which is the closure of the divisor on ${\cal M}_{g,1}$ given
by Weierstrass points of smooth curves, is $${g(g+1)\over 2} \omega -\lambda -
\sum_{i=1}^{g-1} {(g-i)(g-i+1)\over 2} \delta_i.$$

\proof See [Cuk].  Alternatively, this can be done by the method of test
curves.\qed

We recall the Pl\"ucker formula, which counts the ramification points of
a $g^r_d$ on a smooth curve.  It asserts that if $C$ is a smooth curve of
genus $g$ and $V$ a $g^r_d$, then $\sum_{p \in C} \beta(V,p) =
(r+1)(d+r(g-1))$.  Here $\beta(V,p) = \sum_0^r a_i(V,p) - i$, where the $a_i$
give the sequence of orders of vanishing of $V$ at $p$.

\prop{\plucker}
The classical Pl\"ucker formula remains valid for reducible curves with no
nondisconnecting nodes.  (Of course, ``$g^r_d$'' is understood to mean ``limit
linear series''.  Also, the ramification conditions at the node imposed by
the definition of ``limit linear series'' are not considered as contributing
to ramification.)

\proof By induction, it suffices to consider curves with exactly two
components.  Applying the Pl\"ucker formula
to each component separately, we get a total of $(r+1)(2d+r(g-2))$ there.
However, the definition of limit linear series imposes exactly $(r+1)(d-r)$
conditions at the node, leaving $(r+1)(d+r(g-1))$ in total, as claimed.

\define Let $r$ and $d$ be positive integers such that
$a = g-(r+1)(g-d+r) \ge -1$, and $Z$ a possible ramification sequence for
a $g^r_d$ which sums to $a$.  We define
the divisor $D_{g,r,d,Z}$ on $\mbar{g}{1}$ to be the divisor of curves and
points $(C,p)$ such that $C$ admits a $g^r_d$, $V$, whose ramification sequence
at the marked point equals or exceeds $Z$.

The $D_{g,r,d}$ are the divisors of Brill-Noether type on $\mbar{g}{1}$.
In order to study the $D_{g,r,d}$, we consider the
following two maps to $\mbar{g}{1}$:

\item{1.} Map $\mbar{0}{g+1}$ to $\mbar{g}{1}$ by attaching a fixed elliptic
curve at each of the first $g$ points.
\item{2.} Map $\mbar{2}{1} - \W$ to $\mbar{g}{1}$ by attaching a fixed
general curve of genus $g-2$ with two marked points.

We claim that the images of these maps are disjoint from $D$, and that
this implies the linear dependence claimed.  For the first of these maps,
it is enough to count ramification points; there are not enough to spare
on the component of genus $0$.  In particular, at each point of attachment
we must have $\beta \ge r$ on this component, making at least $gr$.  Even
if the remaining $(r+1)(d-r)-gr = g-(r+1)(g-d+r) = a$ units of ramification
are all concentrated at one point, that isn't quite enough.  For the
second map, this follows from the extended Brill-Noether theorem, [EH,
Thm.\ 1.1].

In order to deduce the dependence from this, it is necessary to determine the
pullback map on divisors induced by the two maps given above.  In both cases,
most of the work takes care of itself, since these maps fit into commutative
diagrams in which one arrow is the map in question, one arrow is a map for
which the pullback is computed in [EH, proof of Thms.\ 2.1, 3.1],
and the other arrow(s) are easily understood.
Specifically, to study the first map, insert it into a commutative diagram
as follows, in which the vertical projections are simply the forgetful maps
as shown:

$$\widen\matrix{\mbar{0}{g+1}&\mapright{f'}&\mbar{g}{1}\cr
\mapdown{\pi_{g+1}}&&\mapdown{\pi_1}\cr
\mbar{0}{g}&\mapright{f}&\mgbar\cr}$$

This diagram being commutative, the two pullback maps on $\Pic \otimes \Q$
must be equal.  Since the pullback maps on the vertical arrows are as described
in \pullback, and the map on the bottom arrow is given by
Eisenbud and Harris, we can easily determine $f'^*$ of any class on
$\mbar{g}{1}$ pulled back from $\mgbar$.
For example, since $\pi_{g+1}^* \circ f^* \lambda = 0$ (indeed, $f^* \lambda
= 0$), as in [EH], and since $\pi_1^* \lambda = \lambda$, it follows that
$f'^* \lambda = 0$, and likewise for $\delta_0$ in place of $\lambda$.

\define Let $\theta_i \in \Pic \mbar{0}{g+1} = 
\sum_{S \in T} \delta_{0,S}$, where $T$ runs over subsets of $\{1,\dots,g+1\}$
of cardinality $i+1$ that contain $g+1$.  Also, let 
$\epsilon_i = \sum_{\card S = i} \delta_{0;S}.$

It is easy to see that $\theta_i + \theta_{g-i} = \pi_{g+1}^*
\epsilon_i.$  (As
usual, the case $2i = g$ is an exception: then the pullback is just
$\theta_i$.)  Since
$\epsilon_i = f^* \delta_i$ for $i > 1$, and since for $i < g-1$ we have that
$f'^* \delta_i$ is supported on $\theta_i$, it follows immediately that
$f'^* \delta_i = \theta_i$ in this range.

On the other hand, $$\pi_{g+1}^* \circ f^* \delta_1 = -\sum_{i=2}^{g-2}
{i(g-i)\over(g-1)}\theta_i,$$ as follows from the calculation of $f^* \delta_1$
in [EH, Thm.\ 3.1].  This, therefore, is $f'^* (\delta_1 + \delta_{g-1})$, and
so $$f^*(\delta_{g-1}) = -\sum_{i=1}^{g-2} {i(g-i)\over(g-1)}
\theta_i.$$  Finally
it is necessary to compute $f'^* \omega$.  The easiest way to do this is
simply to use the fact that $f'^* \W = 0$.  Using Cukierman's formula
(\cukwei), one translates this into the statement
that $$f'^* \omega = \sum_{i=1}^{g-2} {(g-i)(g-i-1)\over g(g-1)}
\theta_i .$$

Next we consider the second map, to which a similar procedure applies,
complete with a similar commutative diagram:

$$\widen\matrix{\mbar{2}{1}&\mapright{g'}&\mbar{g}{1}\cr
\mapdown{=}&&\mapdown{\pi_1}\cr
\mbar{2}{1}&\mapright{g}&\mgbar\cr}$$
This time, according to [EH, Sect.\ 2], we have 
$$\eqalign{g^* \delta_0 &= \delta_0,\cr
g^* \delta_1 &= \delta_1,\cr
g^* \delta_2 &= -\omega,\cr
g^* \lambda &= \lambda = \delta_0/10 + \delta_1/5,\cr
g^* \delta_i &= 0\ {\rm for}\ i > 2.}$$
Therefore we will get 
$$\eqalign{g'^* \delta_0  &= \delta_0,\cr
g'^* \lambda &= \lambda = \delta_0/10 + \delta_1/5,\cr
g^* \delta_{g-1} &= \delta_1,\cr
g^* \delta_{g-2} &= \omega,}$$
with all $\delta$ not yet mentioned going to $0$.  After all, $\pi_1{}_*
\delta_i = \delta_i + \delta_{g-i}$, and as $g^* \delta_i$ reflects the
nodes of the curve of genus $2$, the marked points are not on the genus-$i$
side.  Note also that the relation $\lambda = \delta_0/10 + \delta_1/5$ on
$\mgb{2}$ pulls back to $\mbar{2}{1}$ without any change in its appearance.

To complete the description, we show that $g'^* \omega = 0$.  Using the
push-pull formula, we see that it is enough to show that $\omega \cdot
g'_* C = 0$ for any curve $C$ contained in $\mbar{2}{1}$.  But this is
easy, as $\omega$ is the self-intersection of a constant section on the
component of genus $g-2$.  The variation of the other component, which is
attached at some other point, has no effect on this.

To prove the theorem, we show that the intersection of the kernels has
dimension $2$.  This is most easily seen as follows: if we know the
coefficients of $\delta_0$ and $\lambda$ in a divisor contained in
$\ker g'^*$, that determines its coefficients of $\delta_{g-1}$ and
$\delta_{g-2}$.  However, a knowledge of these two coefficients determines
the coefficient of $\omega$ (in order that the coefficient of $\theta_{g-2}$
in the pullback by $f'$ be $0$), and that forces the coefficients of all of
the rest.  Assuming that the $\theta_i$ are independent, we conclude that
the dimension is at most 2; it is at least 2 because the
Weierstrass and Brill-Noether classes are contained in the kernel.  (Of course,
we have only shown that the Brill-Noether class is in the kernel when it is
the class of an effective divisor, but one can check with no difficulty that
it never survives the map $f'^*$.)

Now we prove that the $\theta_i$ are actually independent in $\Pic
(\mbar{0}{g+1}) \otimes \Q$.  This is not quite trivial: there are relations
between the different $\delta_{0;S}$ on $\mbar{0}{n}$.  Since $\mbar
{0}{4} \isisom \P^1$, for example, and the Picard group of $\P^1$ has rank $1$,
we must have $\delta_{0;\{1,2\}} = \delta_{0;\{1,3\}} = \delta_{0;\{1,4\}}$.
But it is not difficult either.

To start with, only $\theta_1$ has nonzero degree on a fiber of $\pi_{g+1}$,
so its coefficient must be $0$ in any relation.  Then $\theta_{g-1}$ is the
only other $\theta$ with nonzero degree on a fiber of $\pi_1$, so its
coefficient must be $0$ as well.  For the rest, put a fixed number $i$, from
$1$ to $g-3$, of the first $g$ points on a $\P^1$ together with $p_{g+1}$,
attach another $\P^1$ at another point with the remaining marked points, and
consider a curve in which one of these points moves along its component.  This
curve will intersect $\theta_{g-1}$ (where the moving point meets another
marked point on its component), $\theta_{i}$ (generically), and $\theta_{i+1}$
(where the moving point reaches the point of attachment), and the coefficient
of $\theta_{i+1}$ will be $1$.  Therefore, by easy induction on $i$, all the
coefficients of $\theta_i$ from $2$ to $g-2$ in our putative relation
must be $0$, and we are done.

Moreover, it is easy to see that the space spanned by these divisors is
actually of dimension $2$ (not $1$).  If $g+1$ is not prime, this is clear,
because $\W$ has a nonzero coefficient of $\omega$ while the pullback of the
Brill-Noether class from $\mgbar$ does not.  If $g+1$ is prime, it is not
much harder.  (Details will appear in a paper presenting the results of [L].)

\section{The General Case}
\noindent Recall the definition of the divisor class $BN$ on $\mgbar$ as
$$(g+3) \lambda - {g+1 \over 6} \delta_0 - \sum i(g-i) \delta_i.$$  If $g+1$ is
composite, $BN$ is effective, but otherwise it may not be.  We may define
divisors on $\mgnbar$ for any $g, n$ by taking any $r$ and $d$ for which the
expected dimension of the space of $g^r_d$'s on a curve of genus $g$ is $k
(\ge -1)$, and imposing exactly enough ramification conditions separately
at the $n$ distinct points to reduce the expected dimension to $-1$.  This
produces a locus in $\mgnbar$ which may in general have components
of various dimensions.  We will refer to any codimension-1 component of any of
these loci as a {\sl Brill-Noether type} divisor, and to the subspace
of $\Pic \mgnbar \otimes \Q$ generated by all
Brill-Noether type divisors as the {\sl Brill-Noether subspace}. 
Observe that if we pull back
a divisor of Brill-Noether type from $\mbar{g}{n-1}$ to $\mgnbar$, we get
another divisor of Brill-Noether type: the same conditions are imposed at
$p_1, \dots, p_{n-1}$, and none at $p_n$.  The main purpose of this paper
is to study
the classes of these divisors: this is what we will do in this section.

The proofs of the earlier results on linear dependence suggest a way to
proceed, and we will follow it.  As before, we consider maps
$$\mbar{0}{g+n} \mapsto \mgnbar$$ given by attaching a fixed elliptic curve
at each of the first $g$ points, and $$\mbar{2}{1} \mapsto \mgnbar$$ given
by attaching a fixed $n+1$-pointed curve of genus $g$
at a marked point.  And as before, the Pl\"ucker
formula shows that the first map misses all of these divisors, while
the extended Brill-Noether theorem proves that the image of $\mbar{2}{1} -
\W$ by second map has trivial intersection with them.  Therefore, their classes
must lie in the intersections of the kernels of these two maps.  These
facts are not quite enough to characterize the Brill-Noether subspace, though.
We will prove the following theorem instead:

\thm{\bndim2}
For any $g > 2$, the Brill-Noether subspace of $\Pic (\mgnbar) \otimes
\Q$ has dimension $1 + n + {n \choose 2}$, unless $g+1$ is prime and $n = 0$.
In particular, its projection
with respect to the standard basis onto the subspace spanned by $\lambda$,
the $\omega$, and the $\delta_{0;\{i,j\}}$ is an isomorphism.

\proof To start the argument, observe that $BN$ has a nonzero coefficient of
$\lambda$, and the pullback of $\W$ by the map which forgets all but the
$i$\/th point has nonzero coefficient of $\omega_i$ and $\lambda$, while all
of their other coefficients in the space we have projected to are $0$.
This, together with the results of the previous section,
takes care of the cases $n = 0, 1$.

We start by proving that the projection is surjective.
The case $n = 2$ is disposed of as soon as we find a
divisor of Brill-Noether type with nonzero coefficient of $\delta_{0;\{1,2\}}$.
To do this, we need a bit of notation.

\define Let $A(g,m,n) = A'(g,d,n)$ be the number of $g^1_d$'s on a general
curve of
genus $g$ ramified to order $m-1$ at a specified point and to order $n-1$
at an unspecified point, where $2d = g+m+n-1$.

\thm{\avalue} ([L], thm. 3.2)
$$A'(g,d,n) = g! (n^2-1)
\sum_{j=\max(0,m+n-d-1)}^{\min(m-1,n-1,d)} {(m+n-2j-1)\over
((d-m-n+j+1)! (d-j)!)}.$$

(It seems almost certain that this was known long before [L], but I know of no
reference.)

In odd genus $g$, we consider the divisor $D$ on $\mbar{g}{2}$ of curves and
points such that there is a $g^1_{(g+3)/2}$ ramified at both of the points.
On the one hand, there are $c_{(g+1)/2} \ g^1_{(g+3)/2}$'s ramified at a given
point, and each of them is ramified at exactly $3g$ other points, so the degree
on a fiber is $3g c_{(g+1)/2}$.  (Here $c_n$ is the $n$\/th
Catalan number $(2n)!/n!(n+1)!$.)  On the other hand, consider $\pi_1{}_*
D \cdot \delta_{0;\{1,2\}}$.  An easy argument with linear series shows this
to consist of the sum of the pullback of a Brill-Noether divisor from $\mgbar$
and a divisor of Brill-Noether type on $\mbar{g}{1}$, to wit, that of curves
and points where the curve has a $g^1_{(g+3)/2}$ doubly ramified at the marked
point.  The degree of this on a fiber is $A(g,1,3) = 24 {g \choose (g+3)/2}$,
and the coefficient of $\delta_{0;\{1,2\}}$ will be $0$ iff this is equal to
$6g c_{(g+1)/2}$.  Expanding out both sides and multiplying through by
$${((g+1)/2)! ((g+3)/2)! \over 6g!},$$ we get $g(g+1)=(g-1)(g+1)$, which holds
for no positive integer, so the coefficient is nonzero as desired.

In even genus, consider the divisor $D$ on $\mbar{g}{2}$ of curves and points
such that there is a $g^1_{(g+4)/2}$ ramified doubly at the first point and
singly at the second.  When we cut and push forward, we get the divisor
$D_{g;2} + D_{g;4}$ on $\mbar{g}{1}$.  Again, to show that $D$ has nonzero
coefficient of $\delta_{0;\{1,2\}}$, we must prove that the sum of the degrees
on the fibers of $D$ is not equal to the degree on the fibers of its
pushforward; in other words, that $A(g,2,3)+A(g,3,2) \ne A(g,1,2)+A(g,1,4)$.
Writing these out in terms of the formula for $A$ given in \avalue\ and
dividing through by $g!$, we get

$$11 \sum_{j=0}^1 {(4-2j)\over (g/2+j-2)!(g/2-j+2)!} \ne {6\over 
(g/2 -1)! (g/2 + 1)!} + {60 \over(g/2 -2)!(g/2+2)!},$$
which on multiplying through by $(g/2-1)!(g/2+2)!$ becomes $$44(g/2-1)+
22(g/2+2) \ne 6(g/2+2)+60(g/2-1),$$ a statement that is always true.

Then, on $\mgnbar$, we can prescribe the coefficients of the
$\delta_{0;\{i,j\}}$ by pulling back these divisors from $\mbar{g}{2}$ in
appropriate ways.  The $\omega$'s are dealt with by pulling back $\W$ from
$\mbar{g}{1}$, and $\lambda$ by pulling back $BN$ from $\mgbar$.  We must
show, now, that a divisor in the Brill-Noether subspace whose coefficients
of $\lambda$, the $\omega$, and the $\delta_{0;\{i,j\}}$ are all $0$ is $0$.

\lemma{\bntobn} Let $S = \{i,j\} \subset \{1, \dots, n+1\}$.  The map $$\Pic
\mbar{g}{n+1} \otimes \Q \mapsto \Pic \mgnbar \otimes \Q$$ given by $D \mapsto
\pi_j{}_* (D \cdot \delta_{0;S})$ maps the Brill-Noether subspace to the
Brill-Noether subspace.  (Note that this is the pullback map on Picard groups
induced by the map $\mgnbar \mapsto \mbar{g}{n+1}$ which takes the point
representing $(C,p_1,\dots,p_i,\dots,p_n)$ to $(C',p_1,\dots,P_1,\dots,P_2,
\dots,p_n)$, where $C'$ is $C$ with a copy of $\P^1$ attached at $p_i$ and
the $P_i$ are points on this copy of $\P^1$.

\proof It suffices to prove this on a set of generators for the Brill-Noether
subspace.  So let $D$ be a divisor of Brill-Noether type.  Then for a point
of $\delta_{0;S}$ to be contained in $D$ means that there is a limit linear
series on the reducible curve the point corresponds to that satisfies the
necessary ramification conditions.  On the $\P^1$ containing $p_i$ and
$p_j$, there is no choice in the matter: we know what the total order of
vanishing must be at the point of attachment.  This forces the total order
of vanishing on the genus-$g$ component, and it is easily checked
that this results in a codimension-$1$ condition for each way to distribute
the vanishing between base points and ramification.  Thus $D$ maps to a sum
of divisors of Brill-Noether type, and the lemma is proved.\qed

\lemma{\showtriv}
Let $D \in BNS$ be a divisor whose coefficients
of $\lambda$, the $\omega$, and the $\delta_{0;\{i,j\}}$ are all $0$, with
respect to the standard basis.  Then
for any $S$, the coefficient of $\delta_{0;S}$ in $D$ is $0$ as well.  (Note
that we do not yet assert that its coefficient of $\delta_0$ must be $0$).

\proof Let $S \subseteq \{1, \dots, n\}\ (\card S > 2)$, and fix a 
$\card S$-pointed $\P^1$ and a general $(n+1-\card S)$-pointed curve of genus
$g$. Consider a family of curves whose base is isomorphic to $\P^1$, and whose
fiber at a point $P$ has the $\card S$-pointed $\P^1$ attached at $P$ to the
curve of genus $g$ at the first point.  For each element $x \in S$, this
family meets $\delta_{0;S - \{x\}}$ once, and it meets $\delta_{0;S}$ with
multiplicity $2 - \card S$---the section on the curve of genus $g$ is constant,
so has self-intersection $0$, while the section on the $\P^1$ is a diagonal
on $\P^1 \times \P^1$, blown up at $\card S$ points.  On the other hand,
it is plain that
this curve does not meet any of the other boundary components or $\lambda$.
In addition, I claim that the intersection with each $\omega_i$ is $0$.

For $i \notin S$, this is obvious.  For $i \in S$, the self-intersection of
the $i$\/th section is $-1$, so this contributes $1$.  For every element $x$ of
$S$ other than $i$, we get a contribution to $\omega_i$ of $1$ from
$\delta_{0;S-\{x\}}$, thus $\card S - 1$.  Finally, we add the intersection
with $\delta_{0;S}$, which is $2 - \card S$; total $0$.

By the extended Brill-Noether theorem [EH, Thm.\ 1.1], this curve misses all
divisors of Brill-Noether type entirely, and so for every $S$ of cardinality
greater than $2$ we get a relation $$(2 - \card S) d_{0;S} + \sum_{x \in S}
d_{0;S \setminus \{x\}} = 0.$$  (Roman letters are the coefficients
of divisors named by similar Greek letters.)  Therefore, if $D$ is a divisor in
the Brill-Noether subspace with coefficients of $\delta_{0;\{i,j\}}$ equal to
$0$, an induction on $\card S$ proves the lemma.\qed

Now we are ready to settle the case $n = 2$ in \bndim2.  Suppose that $D \in
BNS$ has
all coefficients of $\lambda, \omega,$ and $\delta_{0;\{i,j\}}$ equal to 0,
but suppose that the coefficient of $\delta_{i,S}$ is nonzero.  If $S$ is
empty, then consider $\pi_* \delta_{0;\{1,2\}}\cdot D$---it is in the Brill-Noether
subspace and has a nonzero coefficient of $\delta_{i,\emptyset}$,
contradiction.  We may thus assume that all of the $\delta_{i;\emptyset}$ are
$0$.

We now know that $D$ pulls back to $0$ on $\Pic \mbar{0}{g+2}$.  Define
$\theta_{i;S}$ on $\mbar{0}{g+2}$ to be the sum 
$$\sum_{\scriptstyle T \in \{1,2,\dots,g\} \atop \scriptstyle
\card T = i} \delta_{0;T\cup (S+n)},$$
where $S+n$ is the set obtained by adding $n$ to each element of $S$, so that
$\theta_{i;S}$ is the pullback
of $\delta_{i;S}$.  Let $$D = \sum a_i \delta_{i;\{1\}} + d \delta_0;$$
then $\sum a_i \theta_{i;\{1\}} = 0$.  We show that the $\theta_{i;\{1\}}$ are
linearly independent.

To start, only $\theta_{1;\{1\}}$ has nonzero degree on a fiber of $\pi_{g+1}$,
so it cannot appear in a relation.  Next, only $\theta_{g-1;\{1\}}$ of the
remaining $\theta$ has nonzero degree on a fiber of any of the other projection
maps.  For the rest, fix $1 < i < g-2$, put $i$ of the first $g$ points and
$p_{g+2}$ on a $\P^1$ and attach this $\P^1$ to another $\P^1$ which has
the other $g-i$ of the first $g$ marked points and $p_{g+1}$.  Then let
one of the points vary on the first $\P^1$.  This meets only
$\theta_{g-i;\{1\}}$ and $\theta_{g-i+1;\{1\}}$ of the $\theta$, the latter
with multiplicity $1$.  Again it is immediate that the $\theta$ are in fact
linearly independent.

This proves that all of the $a_i$ are $0$, so $D$ is just a multiple of
$\delta_0$.  That means that $D = 0$, though, because the map $\mbar{2}{1}
\mapsto \mbar{g}{2}$ does not pull any nonzero multiple of $\delta_0$
back to a multiple of $\W$. This completes the proof in the case $n = 2$.

Finally, for general $n$, given $D$ we start by removing its coefficients
of $\lambda$, the $\omega$, and the $\delta_{0;\{i,j\}}$.  Its coefficients
of $\delta_{0;S}$ are automatically $0$; suppose that its coefficient of
$\delta_{i;S}$ is nonzero.  Choose $j, k$ so that either $j,k \in S$ or
$j, k \notin S$, cut with $\delta_{0;\{j,k\}}$, and forget the $k$\/th point.
This produces a smaller $n$ and $D$ with such a nonzero coefficient, and
proceeding in this way we get down to the case $n = 2$, contradiction.
So all the $\delta_{i;S}$ have a coefficient of $0$.  We finish the proof
of the theorem by concluding as before that $\delta_0$ cannot appear.\qed

\frenchspacing
\secnum=-1
\section{References}
\leftskip=20.0pt\parindent=-20.0pt\parskip=10pt
[AC] E. Arbarello, M. Cornalba, {\sl The Picard groups of the moduli spaces
of curves.} Topology {\bf 26}, 153--171, 1987.

[Cuk] F. Cukierman, {\sl Families of Weierstrass points.} Duke Math. J. {\bf
58}, 317--346, 1989.

[EH] D. Eisenbud, J. Harris, {\sl The Kodaira dimension of the moduli space
of curves of genus $\ge 23$.} Inv. Math. {\bf 90}, 359--387, 1987.

[HM] J. Harris, I. Morrison, {\sl Moduli of Curves.} GTM 167.  Springer-Verlag,
1998.

[K] F. Knudsen, {\sl The projectivity of the moduli space of stable curves,
II.  The stacks $\mgn$}, Math. Scand. {\bf 52} no. 2, 161--199, 1983.

[L] A. Logan, {\sl Moduli spaces of curves with marked points.}  Harvard
University doctoral dissertation, June 1999.

\bye